\RequirePackage{ifpdf}
\ifpdf 
\documentclass[pdftex]{sigma}
\else
\documentclass{sigma}
\fi

\begin{document}

\allowdisplaybreaks

\renewcommand{\PaperNumber}{106}

\FirstPageHeading

\ShortArticleName{Classical and Quantum Dynamics on Orbifolds}

\ArticleName{Classical and Quantum Dynamics on Orbifolds}

\Author{Yuri A.~KORDYUKOV}

\AuthorNameForHeading{Yu.A.~Kordyukov}

\Address{Institute of Mathematics,
Russian Academy of Sciences,\\ 112 Chernyshevsky Str., Ufa 450008, Russia}
\Email{\href{mailto:yurikor@matem.anrb.ru}{yurikor@matem.anrb.ru}}

\ArticleDates{Received September 04, 2011, in f\/inal form November 20, 2011;  Published online November 23, 2011}

\Abstract{We present two versions of the Egorov theorem for
orbifolds. The f\/irst one is a straightforward extension of the
classical theorem for smooth manifolds. The second one considers an
orbifold as a singular manifold, the orbit space of a Lie group
action, and deals with the corresponding objects in noncommutative
geometry.}


\Classification{58J40; 58J42; 58B34}

\section{Introduction}

The Egorov theorem is a fundamental fact in microlocal analysis and
mathematical physics. It relates the evolution of pseudodif\/ferential
operators on a compact manifold (quantum obser\-vables) determined by
a f\/irst-order elliptic operator with the corresponding evolution of
classical observables~-- the bicharacteristic f\/low on the space of symbols. This theorem
is the rigorous version of the classical-quantum correspondence in
quantum mechanics.

Let $P$ be an elliptic, f\/irst-order pseudodif\/ferential operator on a
compact manifold $X$ with real principal symbol $p\in S^1(T^*X)$.
Let $f_t$ be the bicharacteristic f\/low of the operator $P$, that is,
the Hamiltonian f\/low of $p$ on $T^*X$. The Egorov theorem states
that, for any pseudodif\/ferential operator $A$ of order $0$ with the
principal symbol $a\in S^0(T^*X)$, the operator
$A(t)=e^{itP}Ae^{-itP}$ is a pseudodif\/ferential operator of order
$0$. The principal symbol $a_t\in S^0(T^*X)$ of this operator is
given by the formula
\[
a_t(x,\xi)=a(f_t(x,\xi)), \qquad (x,\xi)\in T^*X\setminus \{0\}.
\]
Here $S^m(T^*X)$ denote the space of smooth functions on
$T^*X\setminus\{0\}$, homogeneous of degree $m$ with respect to a
f\/iberwise $\mathbb R$-action on $T^*X$.

The main purpose of this paper is to extend the Egorov theorem to
orbifolds. We present two versions of the Egorov theorem. The f\/irst
one is a straightforward extension of the classical theorem. The
second one concerns with noncommutative geometry. It considers an
orbifold as a singular manifold, the orbit space of a Lie group
action, and deals with the corresponding noncommutative objects.

Spectral theory of elliptic operators on orbifolds has received much
attention recently (see, for instance, a brief survey in the
introduction of~\cite{Dryden-G-G-W}). In~\cite{Stanhope-Uribe}, the
Duistermaat--Guillemin trace formula was extended to compact
Riemannian orbifolds. This formula has been applied in~\cite{Guilemin-Uribe-Wang} to an inverse spectral problem on some
orbifolds. We believe that our results will play an important role
for the study of spectral asymptotics for elliptic operators on
orbifolds, in particular, for the study of problems related to
quantum ergodicity.

\section{Classical theory}
\subsection{Orbifolds}
In this section we brief\/ly review the basic def\/initions concerning
orbifolds. For more details on orbifold theory we refer the reader
to the book~\cite{Adem-Leida-Ruan}.

Let $X$ be a Hausdorf\/f topological space. An $n$-dimensional
orbifold chart on $X$ is given by a triple $(\tilde{U},G_U,\phi_U)$,
where $\tilde{U}\subset {\mathbb R}^n$ is a connected open subset,
$G_U$ is a f\/inite group acting on~$\tilde{U}$ smoothly and $\phi_U :
\tilde{U} \to X$ is a continuous map, which is $G_U$-invariant
($\phi_U\circ g=\phi_U$ for all $g\in G_U$) and induces a
homeomorphism of $\tilde{U}/G_U$ onto an open subset
$U=\phi_U(\tilde{U})\subset X$. An embedding $\lambda :
(\tilde{U},G_U,\phi_U)\to (\tilde{V},G_V,\phi_V)$ between two such
charts is a smooth embedding $\lambda : \tilde{U}\to \tilde{V}$ with
$\phi_V\circ\lambda=\phi_U$.

An orbifold atlas on $X$ is a family ${\mathcal
U}=\{(\tilde{U},G_U,\phi_U)\}$ of orbifold charts, which cover $X$
and are locally compatible: given any two charts
$(\tilde{U},G_U,\phi_U)$ for $\phi_U(\tilde{U})=U\subset X$ and
$(\tilde{U},G_V,\phi_V)$ for $\phi_V(\tilde{V})=V\subset X$, and a
point $x\in U\cap V$, there exists an open neighborhood $W$ of $x$
and a chart $(\tilde{W},G_W,\phi_W)$ for $W$ such that there are
embeddings $\lambda_U:(\tilde{W},G_W,\phi_W)\hookrightarrow
(\tilde{U},G_U,\phi_U)$ and
$\lambda_V:(\tilde{W},G_W,\phi_W)\hookrightarrow
(\tilde{V},G_V,\phi_V)$.

An (ef\/fective) orbifold $X$ of dimension $n$ is a paracompact
Hausdorf\/f topological space equipped with an equivalence class of
$n$-dimensional orbifold atlases.

Throughout in the paper, $X$ will denote a compact orbifold, $\dim
X=n$.

Let $x\in X$ and $(\tilde{U},G_U,\phi_U)$ be an orbifold chart such
that $x\in U=\phi_U(\tilde{U})$. Take any $\tilde{x}\in \tilde{U}$
such that $\phi_U(\tilde{x})=x$. Let $G_{\tilde{x}}\subset G_U$ be
the isotropy group for $\tilde{x}$. Up to conjugation, this group
doesn't depend on the choice of chart and will be called the local
group at $x$. A point $x\in X$ is called regular if the local group
at $x$ is trivial and it is called singular otherwise. Denote by~$X_{\rm reg}$ the set of regular points of $X$ and by $X_{\rm sing}$ the set
of singular points of $X$. $X_{\rm reg}$ is an open and dense subset of
$X$ whose induced orbifold structure is a manifold structure.

A function $f:X\to {\mathbb C}$ is smooth if\/f for any orbifold chart
$(\tilde{U},G_U,\phi_U)$ the composition $f\left|_U\right.\circ
\phi_U$ is a smooth function on $\tilde{U}$.

The cotangent bundle $T^*X$ of $X$ is an orbifold whose atlas is
constructed as follows. Let $(\tilde{U},G_U,\phi_U)$ is an orbifold
chart over $U\subset X$. Consider the local cotangent bundle
$T^*\tilde{U}=\tilde{U}\times {\mathbb R}^n$. It is equipped with a
natural action of the group $G_U$. The projection map
$T^*\tilde{U}\to \tilde{U}$ is $G_U$-equivariant, so we obtain a
natural projection $p_U:T^*\tilde{U}/G_U\to U$, whose f\/iber
$p_U^{-1}(x)$ is homeomorphic to ${\mathbb R}^n/G_{\tilde{x}}$.
$T^*X$ is obtained by gluing together the bundles $T^*\tilde{U}/G_U$
over each chart $U$ in the atlas of $X$. Namely, let
$(\tilde{U},G_U,\phi_U)$ and $(\tilde{V},G_V,\phi_V)$ be two
orbifold charts for $\phi_U(\tilde{U})=U\subset X$ and
$\phi_V(\tilde{V})=V\subset X$ respectively, and let $x$ belong to
$U\cap V$. There exists an open neighborhood $W$ of $x$ and a chart
$(\tilde{W},G_W,\phi_W)$ for $W$ such that there are embeddings
$\lambda_U:(\tilde{W},G_W,\phi_W)\hookrightarrow
(\tilde{U},G_U,\phi_U)$ and
$\lambda_V:(\tilde{W},G_W,\phi_W)\hookrightarrow
(\tilde{V},G_V,\phi_V)$. These embeddings give rise to
dif\/feomorphisms $\lambda_U: \tilde{W}\to \lambda_U(\tilde{W})\subset
\tilde{U}$ and $\lambda_V: \tilde{W}\to \lambda_V(\tilde{W})\subset
\tilde{V}$, which provide an equivariant dif\/feomorphism
$\lambda_{UV}=\lambda_V\lambda_U^{-1} : \lambda_U(\tilde{W})\to
\lambda_V(\tilde{W})$, the transition function. Then the bundles
$p_U:T^*\tilde{U}/G_U\to U$ and $p_V:T^*\tilde{V}/G_V\to V$ are
glued by the cotangent map $T^*\lambda_{VU} :
T^*\lambda_U(T^*\tilde{W})\to T^*\lambda_V(T^*\tilde{W})$.

\looseness=1
Like in the manifold case, the cotangent bundle $T^*X$ carries a
canonical symplectic structure. Here by a symplectic form on an
orbifold $Y$ we will understand an orbifold atlas ${\mathcal
U}=\{(\tilde{U},G_U,\phi_U)\}$ together with a $G_U$-invariant
symplectic form $\omega_U$ on $\tilde{U}$ for each
$(\tilde{U},G_U,\phi_U)\in {\mathcal U}$ such that, for any
transition function $\lambda_{UV} : \lambda_U(\tilde{W})\subset
\tilde{U}\to \lambda_V(\tilde{W})\subset\tilde{V}$ as above, we have
$\lambda_{UV}^*\omega_V=\omega_U$. An orbifold $Y$ equipped with a
symplectic form $\omega$ is called a symplectic orbi\-fold. The
symplectic structure on $T^*X$ can be constructed as follows.
Consider the orbi\-fold chart $(T^*\tilde{U},G_U,T^*\phi_U)$ induced
by an orbi\-fold chart $(\tilde{U},G_U,\phi_U)$. Then $T^*\tilde{U}$
carries a~canonical symplectic form $\omega_{T^*\tilde{U}}$, which
is invariant with respect to the lifted $G_U$-action. These
symplectic forms are compatible for two dif\/ferent orbifold charts
and def\/ine a symplectic form on~$T^*X$.

The f\/low $F_t$ on a symplectic orbifold $(Y,\omega)$ is Hamiltonian
with a Hamiltonian $H\in C^\infty(Y)$ if, in any orbifold chart
$(\tilde{U},G_U,\phi_U)$, the inf\/initesimal generator $X_H\in
{\mathcal X}(\tilde U)$ of the f\/low satisf\/ies a standard relation
\[
i(X_H)\omega_U=d(H\left|_U\right.\circ \phi_U).
\]
Since $Y$ is not a manifold, this equation can not be reduced to a
system of f\/irst-order ordinary dif\/ferential equations on a manifold.
Nevertheless, one can show the existence and uniqueness of the
Hamiltonian f\/low with an arbitrary Hamiltonian $H$ (for instance,
using quotient presentations, see below). The f\/low $F_t$ leaves the
singular set of the orbifold $Y$ invariant, and its restriction to
the regular part $Y_{\rm reg}$ of $Y$ is the Hamiltonian f\/low of the
function $H\left|_{Y_{\rm reg}}\right.$ in the usual sense. We refer the
reader to \cite{Sjamaar-Lerman} for more information on Hamiltonian
dynamics on singular symplectic spaces.

\subsection{Pseudodif\/ferential operators on orbifolds}
Here we recall basic facts about pseudodif\/ferential operators on
orbifolds (see \cite{Bucicovschi,Girbau-Nicolau1,Girbau-Nicolau2}
for details). We start with some information about orbibundles.

A (real) vector orbibundle over an orbifold $X$ is given by an
orbifold $E$ and a surjective continuous map $p:E\to X$ such that,
for any $x_0 \in X$, there exists an orbifold chart $(\tilde{U},
G_U, \phi_U)$ over $\phi_U(\tilde{U})=U\subset X$ with $x_0 \in U$
and an orbifold chart $(\tilde{U}\times {\mathbb R}^k, G_U,
\tilde{\phi}_U)$ over $\tilde{\phi}_U(\tilde{U}\times {\mathbb
R}^k)=p^{-1}(U)\subset E$ (called a local trivialization of $E$ over
$(\tilde{U}, G_U, \phi_U)$) such that:
\begin{enumerate}\itemsep=0pt
  \item[1)] the action of $G_U$ on $\tilde{U}\times {\mathbb
R}^k$ is an extension of the action of $G_U$ on $\tilde{U}$ given by
\[
g(x,v)=(gx,\rho(x,g)v), \qquad x\in\tilde{U}, \quad v\in{\mathbb R}^k,
\]
where $\rho$ is a smooth map from $\tilde{U}\times G_U$ to the
algebra ${\mathcal L}({\mathbb R}^k)$ of linear maps in ${\mathbb
R}^k$ satisfying
\[
\rho(gx,h)\circ \rho(x,g)=\rho(x,hg),\qquad g,h\in G_U, \quad x\in
\tilde{U};
\]
(in other words, $G_U$ acts by vector bundle isomorphisms of the
trivial vector bundle ${\rm pr}_1: \tilde{U}\times {\mathbb R}^k\to
\tilde{U}$);
  \item[2)] the map ${\rm pr}_1: \tilde{U}\times {\mathbb
R}^k\to \tilde{U}$ satisf\/ies $\phi_U \circ {\rm pr}_1 = p\circ
\tilde{\phi}_U$.
\end{enumerate}
Moreover, any two local trivializations are compatible in a natural
way.

The tangent bundle and the cotangent bundle of an orbifold $X$ are
examples of real vector orbibundles over $X$.

A section $s:X\to E$ is called $C^\infty$, if for each local
trivialization $(\tilde{U}\times {\mathbb R}^k, G_U,
\tilde{\phi}_U)$ over an orbifold chart $(\tilde{U}, G_U, \phi_U)$
the restriction $s\left|_U\right.$ is covered by a $G_U$-invariant
smooth section $\tilde{s}_U: \tilde{U}\to \tilde{U}\times {\mathbb
R}^k$: $s\left|_U\right.\circ \phi_U=\tilde{\phi}_U\circ
\tilde{s}_U$. We denote by~$C^\infty(X,E)$ the space of smooth
section of~$E$ on~$X$.

Now we turn to pseudodif\/ferential operators. Let $X$ be a compact
orbifold, and $E$ a complex vector orbibundle over $X$. A linear
mapping $P : C^\infty(X,E) \to C^\infty(X,E)$ is a (pseudo)
dif\/ferential operator on $X$ of order $m$ if\/f:
\begin{enumerate}\itemsep=0pt
\item[(1)] the Schwartz kernel of $P$ is smooth outside of a neighborhood
of the diagonal in $X\times X$.

\item[(2)] for any $x_0 \in X$ and for any local trivialization
$(\tilde{U}\times {\mathbb C}^k, G_U, \tilde{\phi}_U)$ of $E$ over
an orbifold chart $(\tilde{U}, G_U, \phi_U)$ with $x_0 \in U$, the
operator $C^\infty_c(U,E\left|_U\right.)\ni f \mapsto
P(f)\left|_U\right.\in C^\infty(U,E\left|_U\right.)$ is given by the
restriction to $G_U$-invariant functions of a (pseudo)dif\/ferential
operator $\tilde{P}$ of order $m$ on $C^\infty(\tilde{U},{\mathbb
C}^k)$ that commutes with the induced $G_U$ action on
$C^\infty(\tilde{U},{\mathbb C}^k)$.
\end{enumerate}

All our pseudodif\/ferential operators are assumed to be classical (or
polyhomogeneous), that is, their complete symbols can be represented
as an asymptotic sum of homogeneous components. Denote by
$\Psi^m(X,E)$ the class of pseudodif\/ferential operators of order~$m$
acting on $C^\infty(X,E)$.

It is not hard to show \cite[Proposition 3.3]{Bucicovschi} that the
operator $\tilde{P}$ is unique up to a smoothing operator, so it is
unique if $P$ is a dif\/ferential operator. A pseudodif\/ferential
operator~$\tilde{P}$ on~$\tilde{U}$ that commutes with the action of
$G_U$ has a principal symbol $\sigma(\tilde{P})\in
C^\infty(\tilde{U}\times ({\mathbb R}^n \setminus \{0\}), {\mathcal
L}({\mathbb C}^k))$ that is invariant with respect to the natural
$G_U$-action. One can check that these locally def\/ined functions
determine a global smooth section $\sigma(P)$ of the vector
orbibundle ${\rm End}(\pi^*E)$ on~\mbox{$T^*X \setminus \{0\}$}, the
principal symbol of~$P$. (Here $\pi:T^*X\to X$ is the bundle map and~$\pi^*E$ is the pull-back of the orbibundle $E$ under the map~$\pi$.) The pseudodif\/ferential operator $P$ on $X$ is elliptic if
$\tilde{P}$ is elliptic for all choices of orbifold charts.

\subsection{Classical version of the Egorov theorem}

Let $X$ be a compact orbifold, and $P$ an elliptic, f\/irst-order
pseudodif\/ferential operator on $X$ with real principal symbol $p\in
S^1(T^*X)$. Let $f_t$ be the bicharacteristic f\/low of the operator
$P$, that is, the Hamiltonian f\/low of $p$ on the cotangent bundle
$T^*X$.

As an example, one can consider $P=\sqrt{\Delta_X}$, where
$\Delta_X$ is the Laplace--Beltrami operator associated to a
Riemannian metric $g_X$ on $X$. Its bicharacteristic f\/low is the
geodesic f\/low of the metric $g_X$ on $T^*X$.

The classical version of the Egorov theorem for orbifolds reads as
follows.

\begin{theorem}\label{t:egorov-classical}
For any pseudodifferential operator $A$ of order $0$ with the
principal symbol $a\in S^0(T^*X)$, the operator
\[
A(t)=e^{itP}Ae^{-itP}
\]
is a pseudodifferential operator of order $0$. Moreover, its
principal symbol $a_t\in S^0(T^*X)$ is given~by
\[
a_t(x,\xi)=a(f_t(x,\xi)), \qquad (x,\xi)\in T^*X\setminus 0.
\]
\end{theorem}

The proof of Theorem~\ref{t:egorov-classical} will be given in
Section~\ref{s:proofs}.

\begin{remark}
The classical Egorov theorem plays a crucial role in the proof of
the well-known result due to Shnirelman, stating that the ergodicity
of the bicharacteristic f\/low of a f\/irst-order elliptic
pseudodif\/ferential operator on a compact manifold implies quantum
ergodicity for the operator itself. We will discuss these issues for
orbifolds elsewhere.
\end{remark}

\subsection{The Egorov theorem for matrix-valued operators}

Using the results of \cite{Ja-Str06}, one can easily extend
Theorem~\ref{t:egorov-classical} to pseudodif\/ferential operators
acting on sections of a vector orbibundle over a compact orbifold.

Let $X$ be a compact orbifold, $E$ a complex vector orbibundle on
$X$, and $P$ an elliptic, f\/irst-order pseudodif\/ferential operator
acting on $C^\infty(X,|TX|^{1/2}\otimes E)$ with real scalar
principal symbol $p_1\in S^1(T^*X, {\rm End}(\pi^*E))$,
$p_1(x,\xi)=h(x,\xi){\rm id}_{E_x}$ with $h\in
C^\infty(T^*X\setminus \{0\})$. (Here $|TX|^{1/2}$ denotes the
half-density line orbibundle on $X$.) Let $H_h$ be the associated
Hamiltonian vector f\/ield and $f_t$ the associated Hamiltonian f\/low
on $T^*X$.

Consider a local trivialization $(\tilde{U}\times {\mathbb C}^k,
G_U, \tilde{\phi}_U)$ over an orbifold chart $(\tilde{U}, G_U,
\phi_U)$. Let $\tilde{P}$ be the corresponding $G_U$-invariant,
matrix-valued f\/irst-order pseudodif\/ferential operator on
$\tilde{U}$. The subprincipal symbol of $\tilde{P}$ is a smooth
matrix-valued function ${\rm sub}(\tilde{P})\in
C^\infty(\tilde{U}\times ({\mathbb R}^n\setminus\{0\}),{\mathcal
L}({\mathbb C}^k))$ def\/ined by
\[
{\rm sub}(\tilde{P})(x,\xi)
=p_{0}(x,\xi)-\frac{1}{2i}\sum_{j=1}^n\frac{\partial^2p_1}{\partial
x_j\partial \xi_j}(x,\xi),\qquad x\in \tilde{U}, \quad \xi \in \mathbb
R^n\setminus\{0\},
\]
where $p_k$ is the homogeneous of degree $k$ component in the
asymptotic expansion of the complete symbol of $\tilde{P}$.

If $E$ is trivial, the subprincipal symbol turns out to be well
def\/ined as a function on $T^*X$. In the general case, we consider
the f\/irst-order dif\/ferential operator
\[
\nabla_{H_h}:=H_h+i\,{\rm sub}(\tilde{P})
\]
acting on $C^\infty(\tilde{U}\times ({\mathbb
R}^n\setminus\{0\}),{\mathbb C}^k)$. By \cite{Ja-Str06}, the
operator $\nabla_{H_h}$ is invariantly def\/ined as a covariant
derivative (a partial connection) on the vector orbibundle $\pi^*E$
on $T^*X\setminus \{0\}$ along the Hamiltonian vector f\/ield $H_h$.

This determines a f\/low $\alpha_t$ on $\pi^*E$ by
\[
\alpha_t(x,\xi,v)=(x(t),\xi(t),v(t)), \qquad (x,\xi)\in T^*X\setminus
\{0\}, \quad v\in (\pi^*E)_{(x,\xi)}\cong E_x,
\]
where $(x(t),\xi(t))=f_t(x,\xi)$ and, in local coordinates, $v(t)$
satisf\/ies
\[
\frac{dv(t)}{dt}=i\,{\rm sub}(\tilde{P})(x(t),\xi(t))v(t).
\]
The induced f\/low $\alpha_t^*$ on $C^\infty(T^*X\setminus
\{0\},\pi^*E)$ satisf\/ies
\[
\frac{d}{dt}\alpha_t^*f=\nabla_{H_h}f.
\]
There is also a f\/low $\operatorname{Ad}(\alpha_t)$ on ${\rm
End}(\pi^*E)$, which, in its turn, induces a f\/low
$\operatorname{Ad}(\alpha_t)^*$ on the space $C^\infty(T^*X\setminus
\{0\}, {\rm End}(\pi^*E))$. If $f\in C^\infty(T^*X\setminus \{0\},
{\rm End}(\pi^*E))$,
\[
\frac{d}{dt}\operatorname{Ad}(\alpha_t)^*f=[\nabla_{H_h},f].
\]

\begin{theorem}\label{t:egorov-matrix}
For any $A\in \Psi^0(X,E)$ with the principal symbol $a\in S^0(T^*X,
{\rm End}(\pi^*E))$, the operator
\[
A(t)=e^{itP}Ae^{-itP}
\]
is a pseudodifferential operator of order $0$. Moreover, its
principal symbol $a_t{\in} S^0(T^*X, {\rm End}(\pi^*E))\!$ is given by
\[
a_t=\operatorname{Ad}(\alpha_t)^*a.
\]
\end{theorem}

\subsection{Quotient presentations}
We will need the following well-known fact from orbifold theory due
to Kawasaki \cite{Kawasaki78,Kawasaki81} (see, for instance,
\cite{Bucicovschi,Moerdijk-Mrcun} for a detailed proof).

\begin{proposition}\label{p:quotient}
Let $M$ be a smooth manifold and $K$ a compact Lie group acting on
$M$ with finite isotropy groups. Then the quotient $X = M/K$ $($with
the quotient topology$)$ has a natural orbifold structure. Conversely,
any orbifold is a quotient of this type.
\end{proposition}

Any representation of an orbifold $X$ as the quotient $X \cong M/K$
of an action of a compact Lie group $K$ on a smooth manifold $M$
with f\/inite isotropy groups will be called a quotient presentation
for $X$. There is a classical example of a quotient presentation for
an orbifold $X$ due to Satake. Choose a Riemannian metric on $X$. It
can be shown that the orthonormal frame bundle $M=F(X)$ of the
Riemannian orbifold $X$ is a smooth manifold, the group $K=O(n)$
acts smoothly, ef\/fectively and locally freely on $M$, and $M/K\cong
X$.

\begin{remark}
More generally, one can consider realizations of an orbifold as the
leaf space of a~foliated manifold with all leaves compact and all
holonomy groups f\/inite (a generalized Seifert f\/ibration). The
holonomy groupoid of such a foliation is a proper ef\/fective
groupoid, which provides a characterization of orbifolds in terms of
groupoids.
\end{remark}

Note that if $X \cong M/K$ is a quotient presentation for $X$, then
the pull-back by the natural projection $M \to X$ is an isomorphism
$C^\infty(X)\cong C^\infty(M)^K$ between the smooth functions on $X$
and the $K$-invariant functions on $M$.

There is the following extension of Proposition~\ref{p:quotient}
observed by Kawasaki (see, for instance, \cite{Bucicovschi} for a
detailed proof).

\begin{proposition}\label{p:quotient-E}
Let $\mathcal E$ be a smooth vector bundle over a smooth manifold
$M$ and $K$ a compact Lie group acting on $\mathcal E$ by vector
bundle isomorphisms such that isotropy groups on $M$ are finite.
Then the quotient map $E=\mathcal E/K \to X=M/K$ has a canonical
structure of a vector orbibundle. Conversely, any vector orbibundle
is a quotient of this type.
\end{proposition}

Moreover (see, for instance, \cite[Proposition 2.4]{Bucicovschi}),
if $X \cong M/K$ is a quotient presentation for $X$, $E$ is a vector
orbibundle on $X$, and $\mathcal E$ is the smooth vector bundle over
$M$ given by Proposition~\ref{p:quotient-E}, then
\[
C^\infty(M,\mathcal E)^K\cong C^\infty(X,E).
\]

\subsection{Quotient presentations of the cotangent bundle}

A quotient presentation $X \cong M/K$ for the orbifold~$X$ gives
rise to a quotient presentation for the cotangent bundle~$T^*X$ of
$X$ in the following way. The action of $K$ on $M$ induces an action
of $K$ on the cotangent bundle~$T^*M$. Denote by $\mathfrak k$ the
Lie algebra of $K$. For any $v\in \mathfrak k$, denote by~$v_M$ the
corresponding inf\/initesimal generator of the $K$-action on~$M$. For
any $x\in M$, vectors of the form~$v_M(x)$ with $v\in\mathfrak k$
span the tangent space $T_x(Kx)$ to the $K$-orbit, passing through~$x$. Denote
\[
(T^*_KM)_x  =\big\{\xi\in T^*_xM :\langle \xi, v_M(x)\rangle =0 \
\text{for any}\  v\in \mathfrak k\big\}.
\]
Since the action is locally free, the disjoint union
\[
T^*_KM=\bigsqcup_{x\in M}(T^*_KM)_x
\]
is a subbundle of the cotangent bundle $T^*M$, called the conormal
bundle. The conormal bundle~$T^*_KM$ is a $K$-invariant submanifold
of $T^*M$ such that
\[
T^*_KM/K \cong T^*X.
\]
This gives a quotient presentation for $T^*X$.

This construction is a particular case of the symplectic reduction.
Indeed, the $K$-action on~$T^*M$ is a Hamiltonian action with the
corresponding momentum map $J : T^*M\to {\mathfrak k}^*$ given~by
\[
\langle J(x,\xi),v\rangle = \langle \xi, v_M(x)\rangle, \qquad
(x,\xi)\in T^*M, \quad v\in {\mathfrak k}.
\]
Thus, we see that
\[
T^*_KM=J^{-1}(0),
\]
and the quotient $T^*_KM/K$ is the Marsden--Weinstein reduced space
$M_0$ at $0\in {\mathfrak k}^*$ \cite{Marsden-Weinstein}.

Using quotient presentations, one can show the existence of
Hamiltonian f\/lows on $T^*X$. Let $X\cong M/K$ be a quotient
presentation for $X$. Consider a Hamiltonian $H\in C^\infty(T^*X)$
as a~smooth $K$-invariant function on $T^*_KM$. Let $\tilde{H}\in
C^\infty(T^*M)^K$ be an arbitrary extension of~$H$ to a smooth
$K$-invariant function on $T^*M$. Let $\tilde{f}_t$ be the
Hamiltonian f\/low of $\tilde{H}$ on $T^*M$. Since~$\tilde{H}$ is
$K$-invariant, the f\/low $\tilde{f}_t$ preserves the conormal bundle
$T^*_KM$, and its restriction to $T^*_KM$ (denoted also by
$\tilde{f}_t$) commutes with the $K$-action on $T^*_KM$. So the f\/low
$\tilde{f}_t$ on $T^*_KM$ induces a f\/low $f_t$ on the quotient
$T^*_KM/K=T^*X$, which is called the reduced f\/low. One can show that
this f\/low is a Hamiltonian f\/low on $T^*X$ with Hamiltonian $H$.

\subsection{Pseudodif\/ferential operators and quotient presentations}

Let $X$ be a compact orbifold and $E$ a complex vector orbibundle on
$X$. Let $X \cong M/K$ be a quotient presentation for $X$ and let
$\mathcal E$ be the lift of $E$ to a smooth vector bundle over $M$
given by Proposition~\ref{p:quotient-E}. Let us consider
$C^\infty(X,E)$ (resp. $L^2(X,E)$) as a subspace
$C^\infty(M,\mathcal E)^K$ (resp.\ $L^2(M,\mathcal E)^K$) of
$C^\infty(M,\mathcal E)$ (resp.\ $L^2(M,\mathcal E)$), which consists
of $K$-invariant functions on~$M$. Let $\Pi: L^2(M,\mathcal E)\to
L^2(M,\mathcal E)^K\cong L^2(X,E)$ be the orthogonal projection onto
$L^2(M,\mathcal E)^K$ in~$L^2(M,\mathcal E)$. It is clear that
$\Pi(C^\infty(M,\mathcal E))=C^\infty(M,\mathcal E)^K\cong
C^\infty(X,E)$.

For any pseudodif\/ferential operator $B\in \Psi^m(M,\mathcal E)$,
def\/ine its transversal principal symbol $\sigma(B)\in
S^m(T^*_KM,{\rm End}(\tilde\pi^*\mathcal E))$, where
$\tilde\pi:T^*M\setminus\{0\}\to M$ is the bundle map, as the
restriction of the principal symbol of $B$ to $T^*_KM$. If $B$ is
$K$-invariant, then $\sigma(B)$ is $K$-invariant, so it can be
identif\/ied with an element of the space $S^m(T^*X,{\rm
End}(\pi^*E))$.

We have the following fact, relating pseudodif\/ferential operators on
$M$ and $X$ (see \cite[Proposition 3.4]{Bucicovschi} and
\cite[Proposition 2.3]{Stanhope-Uribe}).

\begin{proposition}\label{p:quotientPsi}
Given a linear operator $A:C^\infty(X,E)\to C^\infty(X,E)$, $A$ is a
pseudodifferential operator on $X$ iff there exists a
pseudodifferential operator $\tilde{A}:C^\infty(M,\mathcal E)\to
C^\infty(M,\mathcal E)$ $($of the same order as $A)$ which commutes
with $K$, such that $\Pi \tilde{A}\Pi=A$.

The transverse principal symbol $\sigma(\tilde{A})\in
S^0(T^*_KM,{\rm End}(\tilde\pi^*\mathcal E))^K$ of $\tilde{A}$
coincides with the principal symbol $\sigma(A)\in S^0(T^*X,{\rm
End}(\pi^*E))$ of $A$ under the identification $C^\infty(T^*_KM,{\rm
End}(\tilde\pi^*\mathcal E))^K\cong C^\infty(T^*X,{\rm
End}(\pi^*E))$.
\end{proposition}

Now we are ready to complete the proofs of
Theorems~\ref{t:egorov-classical} and~\ref{t:egorov-matrix}.

\subsection{Proofs of
Theorems~\ref{t:egorov-classical}
and~\ref{t:egorov-matrix}}\label{s:proofs}

Let $P\in \Psi^1(X)$ be an elliptic operator on $X$ with real
principal symbol $p\in S^1(T^*X)$ and $A\in \Psi^0(X)$. Let $X \cong
M/K$ be a quotient presentation for $X$. Take $\tilde{P}\in
\Psi^1(M)$ and $\tilde{A}\in \Psi^0(M)$ as in
Proposition~\ref{p:quotientPsi}. Without loss of generality, we can
also assume that $\tilde{P}$ is elliptic and its principal symbol is
real. So we have
\[
\Pi\tilde{P}\Pi=P, \qquad \Pi\tilde{A}\Pi=A.
\]
Since $\tilde{P}$ is $K$-invariant, we have
$\tilde{P}\Pi=\Pi\tilde{P}$. Using these facts, we easily derive
that
\[
e^{itP}Ae^{-itP}=\Pi e^{it\tilde{P}}\tilde{A}e^{-it\tilde{P}}\Pi.
\]
By the classical Egorov theorem, the operator
$e^{it\tilde{P}}\tilde{A}e^{-it\tilde{P}}$ is in $\Psi^0(M)$. Since
it commutes with the $K$-action, we conclude that
$e^{itP}Ae^{-itP}\in \Psi^0(X)$. Moreover, for the transversal
principal symbol of $e^{it\tilde{P}}\tilde{A}e^{-it\tilde{P}}$, we
have
\[
\sigma\big(e^{it\tilde{P}}\tilde{A}e^{-it\tilde{P}}\big)(\nu)=
\sigma(\tilde{A})(\tilde{f}_t(\nu)),\qquad \nu\in
T^*_KM\setminus\{0\},
\]
where $\tilde{f}_t$ is the Hamiltonian f\/low of $\tilde{p}$, the
principal symbol of $\tilde{P}$.

By Proposition~\ref{p:quotientPsi}, we have
\[
\sigma\big(e^{it\tilde{P}}\tilde{A}e^{-it\tilde{P}}\big)(\nu)=
\sigma\big(e^{itP}Ae^{-itP}\big)(\nu),\qquad \nu\in
T^*_KM/K\setminus\{0\}\cong T^*X\setminus\{0\},
\]
that immediately completes the proof of
Theorem~\ref{t:egorov-classical}.

The proof of Theorem~\ref{t:egorov-matrix} is similar. Under the
assumptions of Theorem~\ref{t:egorov-matrix}, let $X \cong M/K$ be a~quotient presentation for $X$ and $\mathcal E$ is the lift of $E$ to
a $K$-equivariant smooth vector bundle on $M$ given by
Proposition~\ref{p:quotient-E}. Take $\tilde{P}\in \Psi^1(M,\mathcal
E)$ and $\tilde{A}\in \Psi^0(M,\mathcal E)$ as in
Proposition~\ref{p:quotientPsi}. Without loss of generality, we can
assume that $\tilde{P}$ is elliptic and its principal symbol is
scalar and real. Let~$\tilde{f}_t$ be the Hamiltonian f\/low of
$\tilde{h}\in C^\infty(T^*M\setminus\{0\})$, the principal symbol of~$\tilde{P}$. The subprincipal symbol of $\tilde{P}$ is invariantly
def\/ined as a partial connection $\nabla_{H_{\tilde{h}}}$ on the
vector orbibundle $\tilde{\pi}^*\mathcal E$ along the Hamiltonian
vector f\/ield $H_{\tilde h}$. Therefore, we have the f\/low~$\tilde\alpha_t^*$ on~$C^\infty(T^*M\setminus\{0\}, {\rm
End}(\tilde{\pi}^*\mathcal E))$, which satisf\/ies
\[
\frac{d}{dt}\tilde\alpha_t^*f=\nabla_{H_{\tilde h}}f, \qquad f\in
C^\infty(T^*M\setminus\{0\}, {\rm End}(\tilde{\pi}^*\mathcal E)),
\]
and the f\/low $\operatorname{Ad}(\tilde \alpha_t)^*$ on
$C^\infty(T^*M\setminus\{0\}, {\rm End}(\tilde{\pi}^*\mathcal E))$,
which satisf\/ies
\[
\frac{d}{dt}\operatorname{Ad}(\tilde \alpha_t)^*f=[\nabla_{H_{\tilde
h}},f], \qquad f\in C^\infty(T^*M\setminus\{0\}, {\rm
End}(\tilde{\pi}^*\mathcal E)).
\]
By $K$-invariance, the restriction of $\tilde \alpha_t^*$ to
$C^\infty(T^*_KM\setminus\{0\}, \tilde{\pi}^*\mathcal E)$ takes
$C^\infty(T^*_KM\setminus\{0\}, \tilde{\pi}^*\mathcal E)^K$ to
itself, and there is the following commutative diagram:
\[
  \begin{CD}
C^\infty(T^*_KM\setminus\{0\}, \tilde{\pi}^*\mathcal E)^K @>\tilde
\alpha_t^*>> C^\infty(T^*_KM\setminus\{0\},
\tilde{\pi}^*\mathcal E)^K\\
@VV\cong V @VV\cong V
\\ C^\infty(T^*X\setminus\{0\},\pi^*E)@>\alpha_t^*>> C^\infty(T^*X\setminus\{0\},\pi^*E)
  \end{CD}
\]
A similar statement holds for $\operatorname{Ad}(\tilde
\alpha_t)^*$. Taking into account these facts,
Theorem~\ref{t:egorov-matrix} is a direct consequence of Egorov's
theorem in \cite{Ja-Str06}.

\section{Noncommutative geometry}

\subsection{The operator algebras associated with a quotient orbifold}
Let $X$ be a compact orbifold. The choice of quotient presentation
$X\cong M/K$ for $X$ allows us to consider $X$ as the orbit space of
a Lie group action, which is a typical object of noncommutative
geometry. So we can use some notions and ideas of noncommutative
geometry.

First, one can consider the smooth crossed product algebra
$C^\infty(M)\rtimes K$. As a linear space, $C^\infty(M)\rtimes
K=C^\infty(M\times K)$. The product in $C^\infty(M\times K)$ is
given, for any functions $f, g\in C^\infty(M\times K)$, by
\begin{equation}\label{e:alg1}
(f\ast g)(x,k)=\int_Kf(x,h)g\big(h^{-1}x, h^{-1}k\big)\,dh, \qquad (x,k)\in
M\times K,
\end{equation}
the involution is given, for a function $f\in C^\infty(M\times K)$,
by{\samepage
\begin{equation}\label{e:alg2}
f^*(x,k)=\overline{f(k^{-1}x, k^{-1})}, \qquad  (x,k)\in M\times
K.
\end{equation}
Here $dh$ denotes a f\/ixed bi-invariant Haar measure on $K$.}

It is useful to know that the crossed product algebra
$C^\infty(M)\rtimes K$ is associated with a certain groupoid, the
transformation groupoid, $G=M\rtimes K$. As a set, $G=M\times K$. It
is equipped with the source map $s: M\times K\to M$ given by
$s(x,k)=k^{-1}x$ and the target map $r: M\times K\to M$ given by
$r(x,k)=x$.

For any $x\in M$, there is a natural $\ast$-representation of the
algebra $C^{\infty}(M\times K)$ in the Hilbert space $L^2(K,dk)$
given, for $f\in C^{\infty}(M\times K)$ and $\zeta \in L^2(K,dk)$,
by
\[
R_x(f)\zeta(k)=\int_{K}f\big(k^{-1}x,k^{-1}k_1\big)\,\zeta(k_1)\, dk_1.
\]
The completion of the involutive algebra $C^{\infty}(M\times K)$ in
the norm
\[
\|f\|=\sup_{x\in M}\|R_x(f)\|
\]
is called the reduced crossed product $C^*$-algebra and denoted by
$C(M)\rtimes_r K$.

Since $K$ is compact, this algebra coincides with the full crossed
product $C^*$-algebra $C(M)\rtimes K$, which is def\/ined as the
completion of $C^{\infty}(M\times K)$ in the norm
\[
\|k\|_{\text{max}}=\sup \|\pi(k)\|,
\]
where supremum is taken over the set of all $\ast$-representations
$\pi$ of the algebra $C^{\infty}(M\times K)$ in Hilbert spaces.

There is also a natural representation of $C^\infty(M\times K)$ in
$L^2(M)$ def\/ined for $f\in C^\infty(M\times K)$ and $u\in L^2(M)$ by
\begin{equation}\label{e:Rf}
R(f)u(x)=\int_Kf(x,k)u\big(k^{-1}x\big)dk, \qquad x\in M.
\end{equation}
This representation extends to a $\ast$-representation of
$C(M)\rtimes_r K$.

The $C^*$-algebra $C(M)\rtimes_r K$ can be naturally called the
orbifold $C^*$-algebra associated to the quotient presentation
$X\cong M/K$. From the point of view of noncommutative geometry,
this algebra is a noncommutative analogue of the algebra of
continuous functions on the quotient space $M/K$. It is Morita
equivalent to the commutative algebra $C(X)$.

\subsection{Noncommutative pseudodif\/ferential operators on
orbifolds}\label{s:trpdo} Let $X$ be a compact orbifold. One can
introduce the algebra $\Psi^{*}(M/K)$ of noncommutative
pseudodif\/ferential operators on $X$ associated with a quotient
presentation $X \cong M/K$.

First, let us start with a local def\/inition. Constructing an
appropriate slice for the $K$-action on $M$, one can give the
following local description of the quotient map $p:M\to X$ (see, for
instance, \cite[Proposition 2.1]{Bucicovschi} for details).

\begin{proposition}\label{p:slice}
For any $x\in X$, there exists an orbifold chart
$(\tilde{U},G_U,\phi_U)$ defined in a neighborhood $U\subset X$ of
$x$ such that there exists a $K$-equivariant diffeomorphism
\[
p^{-1}(U)\cong K\times_{G_U}\tilde{U}.
\]
\end{proposition}

Recall that, by def\/inition, $K\times_{G_U}\tilde{U}=(K\times
\tilde{U})/G_U$, where $G_U$ acts on $K\times \tilde{U}$ by
\[
\gamma\cdot (k,y)=\big(k\gamma^{-1}, \gamma y\big),\qquad k\in K, \quad y\in
\tilde{U}, \quad \gamma \in G_U,
\]
 and the $K$-action on $K\times_{G_U}\tilde{U}$ is
given by the left translations on $K$.

Now consider an orbifold chart $(\tilde{U},G_U,\phi_U)$ as in
Proposition~\ref{p:slice}. For any $a \in S ^{m} (K\times K\times
\tilde{U}\times {\mathbb R}^{n})$, def\/ine an operator
\[ \bar{A}: \
C^\infty_c(K\times \tilde{U})\to C^\infty(K\times \tilde{U})
\]
by the formula
\begin{equation}\label{loc}
\bar{A}u(k,y)=(2\pi)^{-n} \int e^{i(y-y')\eta}a(k,k',y,\eta)
u(k',y') \,dk^\prime\,dy'\,d\eta,
\end{equation}
where $u \in C^{\infty}_{c}(K\times \tilde{U})$, $k \in K$, $y \in
\tilde{U}$.

Assume that the operator $\bar{A}$ commutes with the action of $G_U$
on $K\times \tilde{U}$. Then $\bar{A}$ def\/ines an operator
\[
A: \ C^\infty_c(K\times_{G_U}\tilde{U})\cong C^\infty_c\big(p^{-1}(U)\big) \to
C^\infty(K\times_{G_U}\tilde{U}) \cong C^\infty\big(p^{-1}(U)\big).
\]
If, in addition, the Schwartz kernel of the operator $\bar{A}$ is
compactly supported in $(K\times \tilde{U})\times (K\times
\tilde{U})$, then the operator $A$ acts from $C^\infty_c(p^{-1}(U))$
to $C^\infty_c(p^{-1}(U))$ and can be extended in a trivial way to
an operator in $C^\infty(M)$. Such an operator will be called an
elementary operator of class $\Psi^{m}(M/K)$ associated to the
orbifold chart $(\tilde{U},G_U,\phi_U)$.

By def\/inition, the  class $\Psi^{m}(M/K)$ consists of all operators
$A$ in $C^{\infty}(M)$, which can be represented in the form
\[
A=\sum_{i=1}^d A_i + K,
\]
where $A_i$ is an elementary operator of class $\Psi^{*}(M/K)$
associated to an orbifold chart $(\tilde{U_i}{,}G_{U_i}{,}\phi_{U_i})\!$
as in Proposition~\ref{p:slice}, $i=1,\ldots,d,$ and $K\in
\Psi^{-\infty}(M)$.

\begin{remark}
When the group $K$ is discrete (and, therefore, f\/inite), the algebra
$\Psi^{*}(M/K)$ is the crossed product algebra $\Psi^{*}(X)\rtimes
K$.
\end{remark}

The principal symbol of the operator $\bar{A}$ given by (\ref{loc})
is a smooth function on $K\times K\times \tilde{U}\times ({\mathbb
R}^n\setminus\{0\})$ given by
\[
\sigma(\bar{A})(k,k^\prime,y,\eta)=a_{m}(k,k^\prime,y,\eta),\qquad
k,k^\prime \in K, \quad y \in \tilde{U}, \quad \eta\in {\mathbb
R}^n\setminus\{0\},
\]
where $a_{m}$ is the degree $m$ homogeneous component of the
complete symbol $a$. This function is $G_U$-invariant and def\/ines a
smooth function on $(K\times K\times \tilde{U}\times ({\mathbb
R}^n\setminus\{0\}))/G_U$. We put
\begin{gather}
\sigma(A)((k,y,\eta),k^\prime)=\sigma(\bar{A})\big(k,
(k^\prime)^{-1}k,y,\eta\big), \nonumber\\
 (k,y,\eta)\in (K\times
\tilde{U}\times ({\mathbb R}^n\setminus\{0\}))/G_U, \quad k^\prime\in K.\label{k-principal}
\end{gather}

Let $S^{m}(T^*_KM\rtimes K)$ be the space of all smooth functions on
$(T^*_KM\setminus \{0\})\rtimes K$ homogeneous of degree $m$ with
respect to the ${\mathbb R}$-action given by the multiplication in
the f\/ibers of the vector bundle $\pi: T^*_KM \rightarrow M$. The
principal symbol $\sigma(A)$ of an operator $A\in \Psi^{m}(M/K)$
given in local coordinates by the formula (\ref{k-principal}) is
globally def\/ined as an element of the space $S^{m}(T^*_KM\rtimes
K)$.

One can introduce the structure of involutive algebra on
$S^{*}(T^*_KM\rtimes K)$ (see \eqref{e:alg1} and \eqref{e:alg2}),
and show that the principal symbol mapping
\[
\sigma: \ \Psi^{m}(M/K)\rightarrow S^{m}(T^*_KM\rtimes K)
\]
satisf\/ies, for any $A\in \Psi^{m_1}(M/K)$ and $B\in
\Psi^{m_2}(M/K)$,
\[
\sigma(AB)=\sigma(A)\sigma(B),\qquad \sigma(A^*)=\sigma(A)^*.
\]

\begin{example}
For any $f\in C^\infty(M\times K)$, the operator $R(f)$ given by
\eqref{e:Rf} belongs to $\Psi^{0}(M/K)$ and its principal symbol
$\sigma(R(f))\in S^{0}(T^*_KM\rtimes K)$ is given by
\[
\sigma(R(f))=\pi^*_Gf,
\]
where $\pi_G : T^*_KM\times K\to M\times K$ is the natural
projection.
\end{example}

\subsection{Noncommutative Egorov theorem}
As above, let $X$ be a compact orbifold and $X \cong M/K$ a quotient
presentation for $X$. Assume that $P$ is an elliptic, f\/irst-order
pseudodif\/ferential operator on $X$ with real principal symbol $p\in
S^1(T^*X)\cong S^1(T^*_KM)^K$. By Proposition~\ref{p:quotientPsi},
there exists an elliptic, f\/irst-order pseudo\-dif\/fe\-ren\-tial operator
$\tilde{P}$ with real principal symbol $\tilde{p}\in S^1(T^*M)^K$,
which commutes with~$K$ and whose restriction to $C^\infty(M)^K$
agrees with~$P$. Denote by $\tilde{f}_t$ the Hamiltonian f\/low of
$\tilde{p}$ on~\mbox{$T^*M\setminus \{0\}$}. Recall that the f\/low
$\tilde{f}_t$ preserves the conormal bundle $T^*_KM$, and its
restriction to~$T^*_KM$ commutes with the $K$-action on~$T^*_KM$.

Def\/ine a f\/low $F_t$ on $T^*_KM\rtimes K$ by the formula
\[
F_t(\nu,k)=(\tilde{f}_t(\nu),k), \qquad (\nu,k)\in T^*_KM\rtimes K.
\]
Observe that, due to $K$-invariance of the f\/low $\tilde{f}_t$, the
induced map $F^*_t$ on $C^{\infty}(T^*_KM\rtimes K)$ is an
involutive algebra automorphism. The one-parameter automorphism
group $F_t^*$ of the algebra $C^{\infty}(T^*_KM\rtimes K)$ is called
the noncommutative bicharacteristic f\/low of the operator $P$
(associated to the quotient presentation $X\cong M/K$).

\begin{remark}
One can show that the f\/low $F_t^*$ is Hamiltonian with respect to
the natural noncommutative Poisson structure on the algebra
$C^{\infty}(T^*_KM\rtimes K)$ and $p$ is the corresponding
Hamiltonian (see \cite{nc-sympl} for more details).
\end{remark}

The noncommutative version of the Egorov theorem for orbifolds reads
as follows.

\begin{theorem}
\label{t:Egorov-noncom} For any $A\in \Psi^{m}(M/K)$ with the
principal symbol $a\in S^{m}(T^*_KM\rtimes K)$, the operator
\[
\Phi_t(A)=e^{it\tilde{P}}Ae^{-it\tilde{P}}
\]
is an operator of class $\Psi^{m}(M/K)$. Moreover, the principal
symbol $a(t)\in S^{m}(T^*_KM\rtimes K)$ of~$\Phi_t(A)$ is given by
\[
a(t)=F^*_t(a).
\]
\end{theorem}

\begin{proof}
The $K$ orbits on $M$ def\/ine a foliation $\mathcal F$ on $M$ with
compact leaves. The holonomy groupoid $G$ of this foliation
coincides with the transformation groupoid $M\rtimes K$. The
conormal bundle of $\mathcal F$ coincides with $T^*_KM$. As it is
well-known in foliation theory, the conormal bundle of the foliation
carries a natural foliation, ${\mathcal F}_N$, on $T^*_KM$, called
the linearized or the lifted foliation. In our case, this foliation
is given by the $K$ orbits on~$T^*_KM$. The holonomy groupoid~$G_{{\mathcal F}_N}$ of~${\mathcal F}_N$ is the transformation
groupoid $T^*_KM\rtimes K$.

The algebra $\Psi^{*}(M/K)$ of noncommutative pseudodif\/ferential
operators associated with the quotient presentation $X \cong M/K$ is
a particular case of the algebra $\Psi^{*,-\infty}(M,{\mathcal F})$
of transversal pseudodif\/ferential operators on the compact foliated
manifold $(M,{\mathcal F})$ introduced in~\cite{noncom}. In this
setting, Theorem~\ref{t:Egorov-noncom} is a straightforward
consequence of the Egorov theorem for transversally elliptic
operators proved in \cite{mpag} (see also \cite{matrix-egorov}).
\end{proof}

\begin{example}
Suppose that $X$ is a compact manifold (considered as an orbifold).
Then a quotient presentation for $X$ is just a $K$-principal bundle
$\phi:M\to X$. Let~$g_X$ be a Riemannian metric on~$X$. Choose a
bi-invariant metric on~$K$ and a connection on the principal bundle
$\phi:M\to X$. There exists a unique $K$-invariant metric on~$X$,
which makes the map $\phi : (M,g_M)\to (X,g_X)$ into a Riemannian
submersion, with the f\/ibers isometric to $K$. Such a metric is
sometimes called the Kaluza--Klein metric of the connection.

Then, for a Hamiltonian $H\in C^\infty(T^*X\setminus \{0\})$ given
by
\[
H(x,\xi)=|\xi |_{g_X},\qquad (x,\xi)\in T^*X,
\]
the Hamiltonian f\/low
$\tilde{f}_t$ on $T^*M$ is the geodesic f\/low of the metric~$g_M$,
and the reduced Hamiltonian f\/low~$f_t$ on~$T^*X$ is the geodesic
f\/low of the metric~$g_X$.

The corresponding quantum dynamics on $L^2(X)$ and $L^2(M)$ are
described by the operators $P=\sqrt{\Delta_X}$ and
$\tilde{P}=\sqrt{\Delta_M}$ respectively, where~$\Delta_X$ and
$\Delta_M$ are the Laplacians of the metrics~$g_X$ and~$g_M$
respectively. It is well known that the operator $\Delta_M$ can be
expressed in terms of Bochner Laplacians acting on sections of
vector bundles over $X$ associated with the principal bundle
$\phi:M\to X$.

In the case when $K=O(n)$ and $M$ is the orthonormal frame bundle
$F(X)$ of~$X$, the restriction of the geodesic f\/low $\tilde{f}_t$ to
$T^*_KM$ is closely related with the frame f\/low on the frame bundle
$F(X)$ on~$X$ (see~\cite{matrix-egorov} for more details).
\end{example}

\begin{remark}
In this case, both classical and quantum dynamical systems are
noncommutative. It would be interesting to extend some basic results
on quantum ergodicity to this setting. For instance, one can
introduce the notion of ergodicity for the bicharacteristic f\/low
$F_t^*$ on the noncommutative algebra $C^{\infty}(T^*_KM\rtimes K)$
and compare this notion with an appropriate notion of quantum
ergodicity for the operator $P$ itself.
\end{remark}

\vspace{-1mm}

\subsection*{Acknowledgements}

The author was partially supported by the Russian Foundation of
Basic Research (grant no. 10-01-00088).

\vspace{-2mm}

\pdfbookmark[1]{References}{ref}
\LastPageEnding

\end{document}